\documentclass[11pt,reqno]{amsart}
\usepackage{amssymb,latexsym}
\usepackage{graphicx}
\usepackage{fancyhdr,amssymb}
\usepackage{url}		

\theoremstyle{plain}
\numberwithin{equation}{section}

\begin{document}
\fancyhead{}
\renewcommand{\headrulewidth}{0pt}
\fancyfoot{}
\fancyfoot[LE,RO]{\medskip \thepage}

\setcounter{page}{1}

\title[Recurrence of product of linear recursive functions]{Recurrence of product of linear recursive functions}
\author{Cheng Lien Lang}
\address{Department Applied of Mathematics\\
                I-Shou University\\
                Kaohsiung, Taiwan\\
                Republic of China}
\email{cllang@isu.edu.tw}
\thanks{}
\author{Mong Lung Lang}
\email{lang2to46@gmail.com}

\begin{abstract}
We study the recurrence of the product of $n$ functions, each of which satisfies the
 recurrence relation $x(m) = A_1 x(m-1) + A_2 x(m-2) + \cdots +A_s x(m-s)$.

\end{abstract}

\maketitle

\vspace{-.8cm}

\section {Introduction}
\noindent
Let $\{x(m) \,:\, m \in\Bbb Z \,\}$ be a function  defined by the following linear recurrence relation

$$x(m) = A_1 x(m-1) + A_2 x(m-2) + \cdots +A_s x(m-s), \eqno(1.1)$$

\medskip
\noindent where $A_1, A_2,\cdots , A_s$ are constants and $A_i \ne 0$ for some $i$. Denoted by $\{ u(m) \,:\, m\in \Bbb Z\,\}$
 the function that satisfies (1.1) and that $u(1)=u(2) =\cdots  = u(s-1)=0$, $u(s)=1.$ The main purpose of this
  article is to determine a recurrence relation for $X(m)$, where $X(m)$ is a product of $n$ functions, each of
   which satisfies (1.1). The closed form (in terms of matrix) of such a recurrence relation for $X(m)$
   can be found in
     Section 2.  This  fact answers a question (and its generalisation) raised by Cooper and Kennedy [CK]
     which was  partially solved by Stinchcombe [S] ({\em Partially solved} is the sense that
      the  characteristic polynomial of his $x_m$ must admit distinct roots). Two  amazing recurrences
       can be found in Section 5 of [S]. The main result of our study can be found in Proposition 2.3.
     Our method is elementary and uses linear algebra only.

\medskip
  Section 3 studies the  recurrence relation of the product of $n$ recursive functions, each of which
   satisfies the recurrence $W_m = pW_{m-1}-qW_{m-2}$. Note
  that such recurrence  has been studied by Jarden [J], Brenann [B] and Cooper and Kennedy [CK].
   The main results can be found in Proposition 3.2 and 3.2.
   Our proof  is direct and is slightly different from [CK].
    Product of two recursive sequences each of which satisfies $W_m = pW_{m-1}-qW_{m-2}$ is the
     most interesting case. The recurrence for such functions can be used to verify various
      identities concerning Fibonacci  and Lucas numbers. See [LL1] and [LL2] for more detail.

 \medskip
     A rather interesting connection (from our point of view) between the identity $L_n^2 -5F_n^2 = 4(-1)^n$ and the
      Galois field of the characteristic polynomial of the recurrence of $X(m) = F_m^n$ can be found
       in Section 4.

\section {The matrix form of the recurrence of $X(m)$}

  \noindent The main purpose of this section is to give a recurrence relation for $X(m)$
   (a product of $n$ functions $x_i(m)$, each of which satisfies (1.1)). We shall start with an example
    which best reveals our idea and strategy.

\medskip
 \noindent
   {\bf  Example 2.1.} Suppose that $x(m)$ satisfies the recurrence $x(m+2)= px(m+1)-qx(m)$.
   Let $X(m) = x(m)^2.$ Then
  $X(m)$  satisfies the following recurrence relation.

   $$X(m+3) = (p^2-q)X(m+2) +(q^2-p^2q)X(m+1)+q^3X(m).\eqno(2.1)$$

\medskip
  \noindent {\em Proof.}
  Denoted by $u(m)$ the recursive function that satisfies $u(m+2)= pu(m+1)-qu(m)$ and $u(0)=0, u(1)=1$.
   One sees easily that

   $$x(m+r+1) = u(m+1)x(r+1)-qu(m)x(r).\eqno(2.2)$$

  \medskip
 \noindent Take the square of the left and right hand side of (2.2) and shift all the squares to the
   left of the equation, one has

  $$x(m+r+1)^2- u(m+1)^2x(r+1)^2-q^2u(m)^2x(r)^2=
  2qu(m+1)u(m)x(r)x(r+1).\eqno(2.3)$$

 \medskip
 \noindent  We note that on the right hand side of  (2.3), the function $A(r)= 2qx(r)x(r+1)$
   is independent of $m$.
  The cases for $m =1$ and 2 are given as follows.

  $$x(r+2)^2- u(2)^2x(r+1)^2-q^2u(1)^2x(r)^2=
  2qu(2)u(1)x(r)x(r+1),\eqno(2.4)$$
 $$x(r+3)^2- u(3)^2x(r+1)^2-q^2u(2)^2x(r)^2=
  2qu(3)u(2)x(r)x(r+1).\eqno(2.5)$$

  \medskip
  \noindent
  We now put these two equations into the following $2\times 2$ matrix.

 $$  Z= \left[\begin{array}{lr}
x(r+2)^2- u(2)^2x(r+1)^2-q^2u(1)^2x(r)^2
 &u(2)u(1) \\
x(r+3)^2- u(3)^2x(r+1)^2-q^2u(2)^2x(r)^2
 &u(3)u(2)\\
\end{array}\right ]  \eqno(2.6)$$

  \medskip
\noindent
 Denoted by $C_i$ the columns of $Z$.
Applying (2.4)and (2.5), the first column is a multiple of the second column.
To be more accurate, $C_1 -2qx(r+1)x(r)C_2 =
C_1 -A(r)C_2 = 0.$
 Hence the determinant
 of the above matrix $Z$ is zero. Note that this holds for all $r$ and that
  $u(1)=1$, $u(2)= p$ and $u(3) = p^2-q$.
  This completes the study of our example.

\medskip

\noindent We shall now study the general case by a simple generalisation of the
 techniques we  presented in  Example 2.1.

 \medskip
 \noindent {\bf Lemma 2.2.} {\em Let $x(m)$ be a function satisfies $(1,1)$. Then
 $$x(m+r+1) = a_1(m) x(r+1)+ a_2(m) x(r) + \cdots +a_s(m) x(r+2-s),\eqno(2.7)$$
 where $a_1(m) = u(m+1), a_i(m) =A_{i} u(m) + A_{i+1}u(m-1) +\cdots + A_su(m-s+i)$ for $i\ge 2.$
  In particular, $a_i(m)$ satisfies $(1.1)$ for every $i$.}

 \medskip

 \noindent {\em Proof.} Apply mathematical induction. \qed

 \medskip

 \noindent   Let $x_1(m), x_2(m), \cdots , x_n(m)$ be functions, each of which satisfies (1.1).
    It follows from Lemma 2.2 that each one of them satisfies (2.7).
    Let $X(m) = x_1(m)x_2(m) \cdots x_n(m)$. Applying (2.7), one has the following.

 $$ X(m+r+1)-\sum_{i=1}^s a_i(m)^nX( r+2-i) =\sum_{\Delta} A(e_1,e_2,\cdots, e_s) \prod _{i}^s a_i(m)^{e_i},\eqno(2.8)
  $$

 \medskip
 \noindent where $\Delta$ is the set of $(e_1, e_2, \cdots, e_s)$'s such that
  $0\le e_i$ for all $i$, $e_1 + e_2 +\cdots +e_s=n$ and $e_i\ne 0$ for at least two $i$'s as we
  have shifted all the $X(r+2-i)$'s to the left of the equation. $\Delta$ has $k$  members (see Remark 2.4).
   $A(e_1, e_2, \cdots, e_s)$ is a function in $x_i(r+j)$'s, where $1\le i\le n$, $2-s\le j \le 1$
   (take $X(m)=x(m)^2$ in Example 2.1 for example, $A(r) = 2qx(r+1)x(r)$).
   It is important to note that $A(e_1, e_2, \cdots , e_s)$ is independent of $m$. For instance,
   in (2.8),
   regardless the value of $m$,
   the coefficient of $a_{s-1}(m)a_{s}(m)^{n-1} $ is given by $x(r+2-(s-1))x(r+2-s)^{n-1}$.
 Identity
    (2.8) for $m =1, 2, 3, \cdots $ take the following forms

{\small $$ X(1+r+1)-\sum_{i=1}^s a_i(1)^nX( r+2-i) =\sum_{\Delta} A(e_1,e_2,\cdots, e_s) \prod _{i}^s a_i(1)^{e_i}\eqno(2.8_1)
  $$

 $$ X(2+r+1)-\sum_{i=1}^s a_i(2)^nX( r+2-i) =\sum_{\Delta} A(e_1,e_2,\cdots, e_s) \prod _{i}^s a_i(2)^{e_i}\eqno(2.8_2)
  $$

 $$ X(3+r+1)-\sum_{i=1}^s a_i(3)^nX( r+2-i) =\sum_{\Delta} A(e_1,e_2,\cdots, e_s) \prod _{i}^s a_i(3)^{e_i}\eqno(2.8_3)
  $$
$$ \vdots$$
}
 \noindent   Set
 $$ c(m) = X(m+r+1)-\sum_{i=1}^s a_i(m)^nX( r+2-i).\eqno(2.9)$$

   \medskip
   \noindent Associate to  $(2.8_{\,m})$ a row vector of length $k+1$ of the following form (the terms $\,\prod _{i}^s a_i(m)^{e_i}$
    are  ordered according to the  lexicographical order of $\Delta$)
$$  v_m= (c(m), a_1(m)^{n-1} a_2(m),\cdots \cdots ,  a_{s-2}(m) a_{s}(m)^{n-1}  , a_{s-1}(m) a_{s}(m)^{n-1}).\eqno(2.10)$$

 \medskip\noindent
  Similar to how we form matrix (2.6) from (2.4) and (2.5),
   we now put these $k+1$ vectors $v_1, v_{1+1}, \cdots, v_{1+k}$ (in this order) into   the following  $(k+1)\times (k+1)$ square matrix which we call it the {\em recurrence matrix} associated to $X(m)$.

\medskip
\noindent
{
$$ Z= \left [
\begin{array}{lcrr}
 c(1)    &\cdots   &  a_{s-2}(1) a_{s}(1)^{n-1}  & a_{s-1}(1) a_{s}(1)^{n-1}  \\
  c(1+1)   &\cdots   &  a_{s-2}(1+1) a_{s}(1+1)^{n-1}  & a_{s-1}(1+1) a_{s}(1+1)^{n-1}  \\
 \hspace{.5cm}\vdots   &\ddots &\vdots \hspace{2cm}& \vdots\hspace{2cm}\\
  c(1+k)   &\cdots   &  a_{s-2}(1+k) a_{s}(1+k)^{n-1}  & a_{s-1}(1+k) a_{s}(1+k)^{n-1}  \\
   \end{array} \right ]\eqno(2.11)
$$}

\medskip
\noindent
{\bf Proposition 2.3.} {\em  Let  $X(m)$ be  a product of $n$ functions, each of
   which satisfies $(1.1)$.
   Then $X(m)$ admits a recurrence relation.
   Let $Y(m)$ be another function which is also a product of $n$ functions, each of which satisfies
    $(1,1)$. Then $X(m)$ and $Y(m)$ satisfy the same recurrence relation.}

 \medskip
 \noindent {\em Proof.} Denoted by $C_1$  the first column of $Z$ and denoted
  by $C_{e_1e_2\cdots e_s}$ the remaining columns.
 Applying identities  $(2.8_1)$-$(2.8_{k+1})$ and  (2.9),
   the first column of the above matrix $Z$ is a linear combination of the remaining
   columns.
   To be more accurate,
    $$C_1- \sum_{\Delta} A(e_1,e_2,\cdots, e_s)C_{e_1e_2\cdots e_s}=0.\eqno(2.12)
  $$

 \noindent    As a consequence, the determinant of $Z$ is zero.
   Hence the
 set of row vectors $\{R_1, R_2, \cdots, R_{k+1}\}$  of $Z$ is a linearly  dependent set.
Hence there exists $\tau_j$ ($0\le j\le k$), not all zero, such that$\sum \tau_iR_i$
 takes the form $(\sum_{j=0}^{k}\tau_jc(1+j), 0,0,\cdots, 0)$.
  Since the fact that the first column is a linear combination of the remaining columns remains valid
   under row operation, we conclude that $\sum \tau_j c(1+j)=0$.
   We note that in (2.7)-(2.9), $r$ is independent of $m$. This implies that
    $0= \sum \tau_j c(1+j)=\sum \tau_j [X(1+r+1)-\sum_{i=1}^s a_i(1)^nX( r+2-i)]$
     for all $r$.
    As a consequence, this gives a recurrence
    relation for $X(m)$.

    \medskip
    \noindent Since $X(m)$ and $Y(m)$ admit the same recurrence matrix $Z$, they satisfy the same
     recurrence relation. \qed

\medskip
\noindent {\bf Remark 2.4.} Apply the Multinomial Theorem for $(y_1+y_2 +\cdots y_s)^n$,
$k = $
{\tiny  $\left (\begin{array} {c} n+s-1\\
n\\
\end{array} \right ) -s$} .

\medskip
\noindent {\bf Remark 2.5.}
(i)
Denoted by {\small $C_1, C_2, \cdots , C_{k+1}$} the column vectors of $Z$.  Let $V= [C_2, C_3, \cdots ,$
$
 C_{k+1}]$.
 In the case the matrix
 $V= [C_2, C_3,\cdots , C_{k+1}]$
  is of rank $k$,
the
   cofactor expansion along the  first column of $Z$ gives a recurrence
    relation of $X(m)$. In the case the rank of $V $ is $r$ and $r \le k-1$,
      one can still
    apply cofactor expansion to  an $(r+1)\times (r+1)$ submatrix of $Z$ to
    get a recurrence relation of $X(m)$.

     \medskip
           \noindent (ii)
        In [S],
    Stinchcombe [S]
    studies the characteristic polynomial of (1.1) and draws the  conclusion about the
     existence of  a recurrence relation for $X(m) = x(m)^n$.
      Note that his method is different from ours where his calculation
     must be carried through over the complex field $\Bbb C$ and the roots of the
      characteristic polynomial of $x(m)$ must be distinct. Two amazing recurrences can be
       found in Section 5 of [S].

\medskip
\noindent 2.1. {\bf Discussion.} (i)
Let $X(m)$ be given as in Proposition 2.3.  Suppose that the recurrence of $X(m)$ involves $t$ terms. Since
     $\{X(1+k),\cdots, X(t+k)\}$ satisfies the same recurrence relation for arbitrary $k$,
      determinant of the following matrix is zero, which can be viewed as a recurrence relation for $X(m)$ as well.
       This generalises a result stated in [W].

   \medskip
$$ \left [
\begin{array}{lccr}
 X(k+1)& X(k+2)   &\cdots     & X(k+t)\\
X(k+t+1)&  X(k+t+2)   & \cdots     & X(k+2t)\\
 \hspace{.3cm} \vdots                 & \vdots  & \ddots     &\vdots \hspace{.7cm}\\
X(k+ t^2-t+1)& X(k+t^2-t+2) &\cdots &X(k+ t^2)\\
   \end{array} \right ] .\eqno(2.13)
$$

\medskip
\noindent
(ii) As one recurrence function may satisfy more than one recurrence relation, one would like to know
whether Proposition 2.3 offers the best recurrence relation for $X(m)$. The answer is {\em No}. Take
 $x(m) = 2x(m-1) -x(m-2) \,(x(0)=x(1) =1)$ for instance, applying Proposition 2.3,
 the recurrence for $ X(m) = x(m)^2$ is $X(m) = 3X(m-1) -3X(m-2) +X(m-2)$. This is obviously not
  the best recurrence for $X(m)$ as $X(m)$ satisfies the recurrence relation $X(m) = 2X(m-1) -X(m-2)$ as well.
  Is this because the characteristic polynomial of $x(m) = 2x(m-1) -x(m-2)$ has repeated root ? Again, the
   answer is {\em No}. Take  $y(m) = 2y(m-1) -y(m-2) \,(y(0)=0, y(1) =1)$ for instance, $Y(m)=y(m)^2$
    satisfies the relation $Y(m) = 3Y(m-1) -3Y(m-2)+ Y(m-2)$ but $Y(m) \ne 2Y(m)-Y(m-1)$. Hence the initial values
     also play some roles.

 \section {Recurrence Relation for powers of $W_m = pW_{m-1}-qW_{m-2}$}

\noindent Let  $a,b, p,q \in \Bbb C$, $q\ne 0$.
 Following the notations of [Ho],
we define the generalised Fibonacci sequence $\{W_m \} = \{ W_m(a,b\, ; p,q)\}$ by
 $W_0=a$, $W_1 = b$,
$$
W_m= pW_{m-1} -qW_{m-2}.\eqno(3.1)$$

 \medskip
\noindent
Obviously the definition can be extended to negative subscripts ;  that is,
 for $n = 1, 2, 3,\cdots, $ define
  $$W_{-m} = (pW_{-m+1} -W_{-m+2})/q.\eqno (3.2)$$

 \medskip
 \noindent
 In the case $a=0$, $b=1$,  we shall denote the sequence $\{ W_m
  (0,1\,;\, p,q)\}$ by
 $\{u_m\}$. Equivalently,
$$u_m =  W_m(0,1\, ; p,q).\eqno(3.3)$$

\medskip
\noindent
Note that
 $u_{-m} = -q^{-m} u_m$ (this simple fact can be proved by mathematical induction) and that
  $$
  W_{m +r+1} = u_{m+1}W_{r+1} -q u_{m}W_{r}
  .\eqno(3.4)$$

\medskip
\noindent
 The purpose of this section is to give a recurrence relation of $X(m)$, a function which is a
  product of $n$ functions, each of which satisfies the recurrence (3.1).
 Applying our results in Proposition 2.3, the recurrence matrix $Z$ (see (2.11)) is an $n\times n$ matrix of the
  following form and det$\,Z=0$ gives us a recurrence for $X(m)$. A detailed study of the recurrence
   det$\,Z=0$ can be found in sections 3.1 and 3.2.

{\small  $$ Z = \left [
\begin{array}{lcccc}
X(r+2) - X(r+1)u_{2}^n- (-q)^nX(r) u_1^n  & u_{2}^{n-1} u_1^{}    &u_{2}^{n-2}u_1^{2}& \cdots     & u_{2}^{}u_1^{n-1}\\
X(r+3) - X(r+1)u_{3}^n- (-q)^nX(r) u_2^n &  u_{3}^{n-1} u_2^{}    &u_{3}^{n-2}u_2^{2}& \cdots     & u_{3}^{}u_2^{n-1}\\
 \hspace{2cm} \vdots                  &           \vdots            & \vdots         & \ddots     & \vdots\\
X(n+r+1) - X(r+1)u_{n+1}^n- (-q)^nX(r) u_n^n      &   u_{n+1}^{n-1} u_{n}^{}&u_{n+1}^{n-2}u_{n}^{2}&\cdots &u_{n+1}^{}u_{n}^{n-1}\\
   \end{array} \right ]\eqno(3.5)
$$}

\medskip
\noindent 3.1 {\bf The nontrivial case.} Throughout this subsection, we shall assume that
$u_1u_2 \cdots u_{n+1} \ne 0$. Denoted by $c_i$, the $i$-th column of the above matrix.
The determinant of $Z$ can be determined as follows.  Let
$R = (u_2^n, u_3^n, \cdots , u_{n+1}^n)^t$,
$S= (u_1^n, u_2^n, \cdots , u_{n}^n)^t$,
$T = (X(r+2), X(r+3), \cdots , $
$X(n+r+1))^t$
 ($t$ for transpose).
  Set
  $ A= (R, c_2, c_3, $
  $\cdots c_n)$, $ B= (S, c_2, c_3, \cdots c_n)$, $C = (T, c_2, c_3, \cdots c_n)$.
Since the determinant
 function is a linear function, one has
$$ \mbox{det}\, Z = \mbox{det} C- X(r+1) \cdot \mbox{det} A  -(-q)^nX(r)\cdot \mbox{det} B = 0.\eqno(3.6)$$

\medskip
\noindent
Applying $(A3)$ of Appendix $A$, Appendix $B$ and $(C2), (C5)$ of Appendix $C$, one has the following.

  \medskip

  \noindent {\bf Lemma 3.1.} {\em Let $n \in \Bbb N$ be fixed and let $A$, $B$ and $C$ be $n \times n$ matrices given as above $($see $(3.5)$
   and $(3.6))$.
   Suppose that $u_1u_2\cdots u_{n+1} \ne 0$.
   Then
 $$\mbox{det}\, C = \sum_{i=1}^n C_i, \,\,\mbox{det}\,B = (-1)^{n-1}\sigma(n-1)\mbox{det}\,
 A_1,\,\,
  \mbox{det}\,A =\sigma(n)\mbox{det}\,
 A_1,\eqno(3.7)$$

 \medskip
 \noindent where $A_1$ is a nonsingular $n \times n$  matrix, $ \sigma (n) = u_{n+1} u_n \cdots u_3u_2 $
  and the quantity $C_i$ is given by the following.

    $$ C_i = \frac{
    X(r+i+1) \sigma(n)\sigma(n-1)\mbox{det}\,A_1\mbox{det}\,Q_{n-1}^i}
  {u_{i+1}u_i (u_{1-i}u_{2-i}\cdots u_{-1}) (u_1 u_2\cdots u_{n-i})
  \mbox{det}\,Q_n^{i}},\eqno(3.8)$$

\medskip

\noindent where  $1\le i\le n$,  $u_{-m} = -q^{-m}u_m$ $(q$ is the fixed constant
 given in $(3.1))$ and $Q_m$ is an $m\times m$ matrix of determinant $ q^{m(m-1)/2}$.}
\medskip

\noindent {\em Proof. } See Appendix $A$, $B$ and $C$.\qed

\medskip

\noindent
Note  that det$\,A_1$ appears in $C_i$, det$\,B$ and det$\,A$ and that det$\,A_1$ is nonzero.
 Applying (3.7), (3.8) and  $u_{-m} = -q^{-m}u_m$, one may now simplify
 (3.6) and conclude that
 $$
  \sum_{i=0}^{n+1}
  (-1)^{i}  q^{i(i-1)/2}
(n+1|i)_u
 X(n+r+1-i) =0,\eqno(3.9)$$
\noindent where
$(m|k)_u =1$
if $k=0$ and $
(m|k)_u =
u_mu_{m-1}\cdots u_{m-k-1}/u_ku_{k-1}\cdots u_1$ if $1\le k \le m$.
Note that $(m|k)_u$ is known as the generalised binomial coefficient.
 In summary, the following proposition holds.

\medskip
\noindent {\bf Proposition 3.2.} {\em  Let $n \in \Bbb N$ be fixed and let $X(m)$ be a product
 of $n$ functions, each of which satisfies $(3.1)$.
 Suppose that $u_1 u_2 \cdots u_{n+1}\ne 0$.
 Then $X(m)$ satisfies the
  following recurrence.
   $$
  \sum_{i=0}^{n+1}
  (-1)^{i} q^{i(i-1)/2} (n+1|i)_u
 X(m-i) =0.\eqno(3.10)$$}

\medskip
\noindent {\em Proof.}  {\bf Alternative Proof of Proposition 3.2.} One sees easily that the $X(m) $ in Proposition 3.2
 and $u_m^n$ admit the same recurrence matrix (see (2.11) for the definition and (3.5) for the actual
  matrix). Hence $X(m)$ and $u_m^n$
  satisfy  the same recurrence relation (see Proposition 2.3). One may now prove the proposition by applying Jarden's Theorem
   which gives the  recurrence for $u_m^n$
   (see Theorem 4.1
   of [CK]).\qed

   \medskip

\noindent 3.2. {\bf The trivial case.} Throughout this section, $u_1 u_2 \cdots u_{n+1} =0$.
 This implies that $u_{k}=0$, for some $k$, where  $1\le k \le n+1$.
Identity (3.4) now becomes $W_{k+r+1} = u_{k+1}W_{r+1}.$
 Since $X(m)$ is a product of $n$ function, each of which satisfies  $W_{k+r+1} = u_{k+1}W_{r+1}$,
  one has
  $X(k+r+1) = u_{k+1}^nX(r+1)$.
  Since this holds  for every  $r \in \Bbb Z$ and $k$ is fixed,
   the following proposition is clear.

 \medskip
\noindent {\bf Proposition 3.3.} {\em  Let $n \in \Bbb N$ be fixed and let $X(m)$ be a product
 of $n$ functions, each of which satisfies $(3.1)$.
 Suppose that $u_k= 0$, for some $k$, where $1\le k\le n+1$.
 Then $X(m)$ satisfies the
  following recurrence.
   $$X(k+r+1)= u_{k+1}^nX(r+1).\eqno(3.11)$$
  }

  \noindent {\bf Example 3.4.} Let $u_0=1, u_1 =1$ and $u_n= 2u_{n-1} -4u_{n-2}$. Then
  $u_2 = 2, u_3 =0, u_4 = -8$. By Proposition 3.3,
   $ u_m^3$
   satisfies the following recurrence

   $$ u_{r+4}^3 + 8^3u_{r+1}^3=0.\eqno(3.12)$$

\medskip
  \noindent Note that $u_n^3$ also satisfies the recurrence relation $u_{r+4}^3 -8u_{r+3}^3 -512u_{r+1}^3 +4096u_r^3=0$.

 \medskip
 \noindent  3.3. {\bf  Application.}
  The following corollary of Proposition 3.2 has been applied by Lang and Lang ([LL1] and [LL2])  to prove various identities concerning
   the generalised Fibonacci sequence.

\medskip
\noindent {\bf Corollary 3.5.} {\em  $ A(m) = W_{2m}$, $B(m) =W_{m}W_{m+r}$ and $C(m) = q^m$  satisfy the following
  recurrence relation}
 {\small  $$X(m+3) = (p^2-q)X(m+2) +(q^2-p^2q)X(m+1)+q^3X(m).\eqno (3.13)$$}

  The  following corollary of Proposition 3.2 shows that the recurrence (3.10) can be used to
 describe  the characteristic polynomial of $Q_v$ (see $(A4)$ of Appendix $A$) which reveals the
  fact that the recurrence relation carries a lot of information  about the generalised
   Fibonacci numbers. Proof of Corollary 3.6 can be found in Appendix $D$.

 \medskip

  \noindent {\bf Corollary 3.6.} {\em Let $Q_{v}$ be given as in $(A4)$ of Appendix $A$.
   Then the characteristic polynomial of $Q_{v}$ is
   $
  \sum_{i=0}^{v}
  (-1)^{i} q^{i(i-1)/2} (v|i)_u x^i.$
  Note that $f(x)$ is the polynomial that characterises the recurrence $(3.10)$.
  }

\section {The Galois Group of the Recurrence Relation}

Let $p$ and $q$ be given as in (3.1) and let

 $$\phi_n(p, q, x)=
  \sum_{i=0}^{n}
  (-1)^{i} q^{i(i-1)/2} (n|i)_u x^i.\eqno(4.1)$$

\noindent We call the polynomial in (4.1) (see Corollary 3.6 also) the Galois polynomial of the recurrence relation (3.10). Similarly we call
 the Galois group $G_n(p, q)$ of $\phi_n(p, q, x)$ over $\Bbb Q$ the Galois group of the recurrence relation (3.10).

 \medskip
 \noindent {\bf Proposition 4.1.} {\em $G_n(1, -1) \cong \Bbb Z_2,$  a cyclic group of order $2$, for all $n \ge 2$.}

 \medskip
 \noindent {\em Proof.} Let $a$ and $b$ be roots of $ \phi_2(x) = x^2-x-1=0$. One  can show by induction
  that (see [B])
  $$\phi_n(1,-1,x) =(ab)^{n-1} (x-a^n)(x-b^n)\phi_{n-2}(1,-1,x/ab)\eqno(4.2)$$

\medskip
  \noindent for all $n\ge 2$ $(\phi_0(x)$ is defined to be 1). Hence $\phi_n(1,-1,x)$ splits completely in
   the field $\Bbb Q
   (\sqrt 5) $ and $G_n(1,-1) \cong \Bbb Z_2$.\qed

 \medskip
 \noindent 4.1. {\bf Discussion.} Let $F_n$ and $L_n$ be the $n$-th Fibonacci and Lucas numbers respectively.
  It is a well known fact that
  $$
 L_n^2 -5F_n^2 = 4(-1)^n .\eqno (4.3)$$

 \noindent Identity (4.3), from our point of view, is a {\em must} rather than {\em an interesting
  connection between Fibonacci and Lucas numbers} as the following suggested.
   Recall another well known fact about the factorisation of $\phi_n(1,-1, x)$.

  $$\phi_n(1,-1,x) = (-1)^n(x^2-L_nx +(-1)^n)\phi_{n-2}(1,-1, -x).\eqno(4.4)$$

 \medskip
  \noindent Since $\phi_n(1,-1,x)$ splits in $\Bbb Q(\sqrt 5)$ and $x^2-L_nx +(-1)^n =0$
   splits in $\Bbb Q(\sqrt{L_n^2 -4(-1)^n})$, one must have
   $$L_n^2 -4(-1)^n = 5A_n^2,\eqno(4.5)$$

   \noindent for some $A_n\in \Bbb N$.
   This tells us that the difference between $4(-1)^n$ and the square of the Lucas number $L_n$ must be five times a square $A_n^2$. As for why $A_n$ must be $F_n$,
    we note that both $F_n^2$ and $L_n^2 $ satisfy the recurrence relation (3.13)
   (to be more accurate, the recurrence $X(n+3) = 2X(n+2)+2X(n+2)-X(n)$) and that
    $L_n^2-4(-1)^n$ and $F_n^2$ have the same initial values. In general, one has
     (Lemma 3.3 of Cooper and Kennedy [CK]),

    $$\phi_{n}(p,q, x) = \prod_{j=0}^{n} (x-\alpha^j\beta^{n-j}),\eqno(4.6)$$

    \noindent where $\alpha$ and $\beta$ are roots of $x^2-px+q =0$. Hence $G(p,q) \cong \Bbb Z_2$
     for all $n \ge 2$. This fact will give an identity similar to (4.5).

 \section{ Appendix $A$}
  Throughout the appendix $u_1 u_2\cdots u_{n+1} \ne0$.
 We shall give full detail of how the determinant of the matrix $A$ is evaluated ($A$ is given as follows).
  Applying our technique given in this appendix, the determinants of $B$ and $C$ can be calculated similarly
   (see Appendix $B$ and $C$).

   \medskip

$$ A= \left [
\begin{array}{ccccc}
u_{2}^n& u_{2}^{n-1} u_1^{}    &u_{2}^{n-2}u_1^{2}& \cdots     & u_{2}^{}u_1^{n-1}\\
u_{3}^n &  u_{3}^{n-1 }u_2  &u_{3}^{n-2}u_2^{2}& \cdots     & u_{3}^{}u_2^{n-1}\\
  \vdots                  &           \vdots            & \vdots         & \ddots     & \vdots\\
u_{n+1}^n &   u_{n+1}^{n-1} u_{n}^{}&u_{n+1}^{n-2}u_{n}^{2}&\cdots &u_{n+1}^{}u_{n}^{n-1}\\
   \end{array} \right ]\eqno(A1)
$$

\bigskip
\noindent Note  that $u_{i+1}$ is a common factor of the entries of the $i$-th row.
 Hence  det$\,A$ can be
 written as the product  $u_{n+1}u_n \cdots u_2\,$det$\,A_1$, where

$$ A_1= \left [
\begin{array}{ccccc}
u_{2}^{n-1}& u_{2}^{n-2} u_1^{}    &u_{2}^{n-3}u_1^{2}& \cdots     & u_1^{n-1}\\
u_{3}^{n-1} &  u_{3}^{n-2 }u_2  &u_{3}^{n-3}u_2^{2}& \cdots     & u_2^{n-1}\\
  \vdots                  &           \vdots            & \vdots         & \ddots     & \vdots\\
u_{n+1}^{n-1} &   u_{n+1}^{n-2} u_{n}^{}&u_{n+1}^{n-3}u_{n}^{2}&\cdots &u_{n}^{n-1}\\
   \end{array} \right ]\eqno(A2)
$$

   \bigskip
   \noindent For our convenience, we shall define $\sigma(n)$, $\tau(n)$   as follows which will be used in the following discussion.
    $$ \sigma (n) = \prod_{i=2}^{n+1} u_i ,\,\,\,  \tau(n) = q^{n(n-1)/2},\,\,\,
    \mbox{det}\,A = \sigma (n)\,\mbox{det}\,A_1.\eqno(A3)$$

\medskip
\noindent
The rest of this section is devoted to the determination of the determinant of $A_1$.
    To save space, we denote the binomial coefficient
{\tiny $
  \left(\begin{array}{c}
n\\
k\\
\end{array}\right )$} by $(n|k)$.
 Consider the following matrix.

\medskip
{\small $$\,\,\,\, Q_n=
 \left [
 \begin{array}{llcrr}
 (n-1|0)p^{n-1}&    (n-2|0)    p^{n-2}& \cdots     &  p  &     1            \\
 (n-1|1) p^{n-2} (-q)& (n-2|1)    p^{n-3}(-q)& \cdots     & -q  &  0    \\
 (n-1|2) p^{n-3} (-q)^2   &(n-2|2)p^{n-4}(-q)^2& \cdots     & 0    & 0        \\
  \hspace{1cm}\vdots                  &      \hspace{1cm}     \vdots            & \ddots         & \vdots     & \vdots\\
  (n-1|n-2)   p (-q)^{n-2}& (n-2|n-2)(-q)^{n-2}&\cdots &0  & 0               \\
  (n-1|n-1) (-q)^{n-1}& 0&\cdots &0 & 0            \\
   \end{array} \right ].
    \eqno(A4)
   $$}

\medskip
\noindent
Note that $Q_n$ is an $n\times n$ matrix where the $i$-th column of $Q_n$  gives the coefficients of the
 binomial expansion of $(px-qy)^{n-i}$ and that det$\,Q_n = \tau(n) = q^{n(n-1)/2}$.
Since $u_{r+1} = p u_r - qu_{r-1}$, $u_0 =0$, $u_1=1$ (see (1.3)),
the multiplication of $A_1$ by $Q_n$  shifts the indices. That is,

\medskip
{\small $$
A_1Q_n^{-1} =
\left [\begin{array}{cccc}
u_{1}^{n-1}& u_{1}^{n-2} u_0^{}    &\cdots     & u_0^{n-1}\\
u_{2}^{n-1} &  u_{2}^{n-2 }u_1  & \cdots     & u_1^{n-1}\\
  \vdots                          & \vdots         & \ddots     & \vdots\\
u_{n}^{n-1} &   u_{n}^{n-2} u_{n-1}^{}&\cdots &u_{n-1}^{n-1}\\
   \end{array} \right ],\,\,
A_1Q_n =
\left [\begin{array}{ccccc}
u_{3}^{n-1}& u_{3}^{n-2} u_2^{}    & \cdots     & u_2^{n-1}\\
u_{4}^{n-1} &  u_{4}^{n-2 }u_3  & \cdots     & u_3^{n-1}\\
  \vdots                  &           \vdots             & \ddots     & \vdots\\
u_{n+2}^{n-1} &   u_{n+2}^{n-2} u_{n+1}^{}&\cdots &u_{n+1}^{n-1}\\
   \end{array} \right ].
  \eqno(A5)
   $$}

\bigskip
\noindent Since $u_0=0, u_1 = 1$, the first row of $A_1Q_n^{-1}$  is $(1^{n-1}, 0,0,\cdots, 0)$.
 Applying the cofactor expansion to the first row of  $(A5)$,  the determinant of $A_1 Q_n^{-1}$ is
 the determinant of the following matrix.

\medskip
$$
Z=
\left [\begin{array}{cccc}
 u_{2}^{n-2} u_1^{}    &u_{2}^{n-3}u_1^{2}& \cdots     & u_1^{n-1}\\
 u_{3}^{n-2 }u_2 &u_{3}^{n-3}u_2^{2}& \cdots     & u_2^{n-1}\\
          \vdots            & \vdots         & \ddots     & \vdots\\
   u_{n}^{n-2} u_{n-1}^{}&u_{n}^{n-3}u_{n-1}^{2}&\cdots &u_{n-1}^{n-1}\\
   \end{array} \right ]
  \eqno(A6)
   $$

\bigskip

\noindent Compare the matrices in $(A1)$ (an $n\times n$ matrix)
 and $(A6)$ (an $(n-1)\times (n-1)$ matrix), we have established a recursive process which enables us
 to calculate
 the determinant of the matrix $A$ as well as $A_1$.
 $$
 \mbox{det}\,A_1 =\tau(n)\tau(n-1)\cdots \tau(2)\sigma(n-2)\sigma(n-3)
 \cdots \sigma(1)\ne 0
 .\eqno(A7)$$

\medskip
\noindent Since the multiplication of $A_1$  by $Q_n$ shifts the indices (see $(A5)$), the following is clear.
Note that this matrix plays an important in the study of matrix $C$ (see Appendix $C$).

\medskip
$$
\hspace{-1cm}
A_1Q_n^{-i}=
\left [
\begin{array}{ccccc}
u_{2-i}^{n-1}& u_{2-i}^{n-2} u_{1-i}^{}    &u_{2-i}^{n-3}u_{1-i}^{2}& \cdots     & u_{1-i}^{n-1}\\
u_{3-i}^{n-1} &  u_{3-i}^{n-2 }u_{2-i}  &u_{3-i}^{n-3}u_{2-i}^{2}& \cdots     & u_{2-i}^{n-1}\\
  \vdots                  &           \vdots            & \vdots         & \ddots     & \vdots\\
u_{1}^{n-1} &   u_{1}^{n-2} u_{0}^{}&u_{1}^{n-3}u_{0}^{2}&\cdots &u_{0}^{n-1}\\
\vdots                  &           \vdots            & \vdots         & \ddots     & \vdots\\
u_{n+1-i}^{n-1} &   u_{n+1-i}^{n-2} u_{n-i}^{}&u_{n+1-i}^{n-3}u_{n-i}^{2}&\cdots &u_{n-i}^{n-1}\\
   \end{array} \right ]\eqno(A8)
$$

 \bigskip
\noindent Note that the $i$-th row of $A_1 Q_n^{-i}$ takes the form $(1^{n-1},0,0, \cdots, 0)$.
 The determinant of $A_1 Q_n^{-i}$ can be determined by the cofactor expansion by the $i$-th row
  (see $(C5)$ of Appendix $C$).

\section {Appendix $B$}

  Throughout the appendix $u_1 u_2\cdots u_{n+1} \ne0$.
The matrix $B$ is given as follows.

$$
B =  \left [
\begin{array}{ccccc}
u_{1}^n& u_{2}^{n-1} u_1^{}    &u_{2}^{n-2}u_1^{2}& \cdots     & u_{2}^{}u_1^{n-1}\\
u_{2}^n &  u_{3}^{n-1 }u_2  &u_{3}^{n-2}u_2^{2}& \cdots     & u_{3}^{}u_2^{n-1}\\
  \vdots                  &           \vdots            & \vdots         & \ddots     & \vdots\\
u_{n}^n &   u_{n+1}^{n-1} u_{n}^{}&u_{n+1}^{n-2}u_{n}^{2}&\cdots &u_{n+1}^{}u_{n}^{n-1}\\
   \end{array} \right ]\eqno(B1)
$$

\medskip
\noindent
After the removal
  of the common factor $u_i$ from each entry of the $i$-th row, one sees that the resulting matrix is
   just the matrix $A_1$ (see equation $(A2)$)
   when the first column is moved to the last.
   Hence det$\,B =(-1)^{n-1} \sigma (n-1)$det$\,A_1$
    (see $A(3)$ for the definition of $\sigma(n-1)$).

\section{Appendix $C$}
  Throughout the appendix $u_1 u_2\cdots u_{n+1} \ne0$.
The matrix $C$ is given as follows. Its determinant can be calculated by the cofactor
 expansion of the first column. The calculation is tedious but elementary.

$$ C= \left [
\begin{array}{ccccc}
W_{r+2}^n& u_{2}^{n-1} u_1^{}    &u_{2}^{n-2}u_1^{2}& \cdots     & u_{2}^{}u_1^{n-1}\\
W_{r+3}^n &  u_{3}^{n-1 }u_2  &u_{3}^{n-2}u_2^{2}& \cdots     & u_{3}^{}u_2^{n-1}\\
  \vdots                  &           \vdots            & \vdots         & \ddots     & \vdots\\
W_{n+r+1}^n &   u_{n+1}^{n-1} u_{n}^{}&u_{n+1}^{n-2}u_{n}^{2}&\cdots &u_{n+1}^{}u_{n}^{n-1}\\
   \end{array} \right ]\eqno(C1)
$$

\bigskip
\noindent
 Denoted by $C_i$ the $i$-th cofactor of the first column of $C$. The determinant of $C$ is given by the
  following.
   $$\mbox{det}\, C= C_1 +C_2 +\cdots +C_n,\eqno(C2)$$

\medskip
\noindent where
  $C_i = (-1)^{i+1}W_{r+i+1}^n$det$\,X$ and $X$ is the  $(n-1)\times (n-1)$ matrix that takes  the
   following form.

{\small $$ X =
\left [
\begin{array}{cccc}
 u_{2}^{n-1} u_1^{}    &u_{2}^{n-2}u_1^{2}& \cdots     & u_2u_1^{n-1}\\
  u_{3}^{n-1 }u_2  &u_{3}^{n-2}u_2^{2}& \cdots     & u_{3}^{}u_2^{n-1}\\
             \vdots            & \vdots         & \ddots     & \vdots\\
   u_{i}^{n-1} u_{i-1}^{}    &u_{i}^{n-2}u_{i-1}^{2}& \cdots     & u_{i}^{}u_{i-1}^{n-1}\\
  u_{i+2}^{n-1 }u_{i+1}  &u_{i+2}^{n-2}u_{i+1}^{2}& \cdots     & u_{i+2}^{}u_{i+1}^{n-1}\\
         \vdots            & \vdots         & \ddots     & \vdots\\
 u_{n+1}^{n-1} u_{n}^{}&u_{n+1}^{n-2}u_{n}^{2}&\cdots &u_{n+1}^{}u_{n}^{n-1}\\
   \end{array} \right ]
   \eqno(C3)
$$}

\bigskip
\noindent Note that $u_{j+1}u_j$ is a common factor of the entries of the $j$-th row when $j\le i-1$ and that
  $u_{j+2}u_{j+1}$ is a common factor of the entries of the $j$-th row when $j \ge i$. As a consequence, the determinant
   of $X$ is the product of $\sigma(n)\sigma(n-1)/u_{i+1}u_i$ and the determinant of $Y$, where

{\small
 $$Y =
 \left [
\begin{array}{cccc}
 u_{2}^{n-2}    &u_{2}^{n-3}u_1^{}& \cdots     & u_1^{n-2}\\
  u_{3}^{n-2 } &u_{3}^{n-3}u_2^{}& \cdots     & u_2^{n-2}\\
             \vdots            & \vdots         & \ddots     & \vdots\\
   u_{i}^{n-2}     &u_{i}^{n-3}u_{i-1}^{}& \cdots     & u_{i-1}^{n-2}\\
  u_{i+2}^{n-2 }  &u_{i+2}^{n-3}u_{i+1}^{}& \cdots     & u_{i+1}^{n-2}\\
         \vdots            & \vdots         & \ddots     & \vdots\\
 u_{n+1}^{n-2} &u_{n+1}^{n-3}u_{n}^{2}&\cdots &u_{n}^{n-2}\\
   \end{array} \right ]\,,\,\,\,
   Y Q_{n-1}^{-i}=
 \left [
\begin{array}{cccc}
 u_{2-i}^{n-2}    &u_{2-i}^{n-3}u_{1-i}^{}& \cdots     & u_{1-i}^{n-2}\\
  u_{3-i}^{n-2 } &u_{3-i}^{n-3}u_{2-i}^{}& \cdots     & u_{2-i}^{n-2}\\
             \vdots            & \vdots         & \ddots     & \vdots\\
   u_{0}^{n-2}     &u_{0}^{n-3}u_{-1}^{}& \cdots     & u_{-1}^{n-2}\\
  u_{2}^{n-2 }  &u_{2}^{n-3}u_{1}^{}& \cdots     & u_{1}^{n-2}\\
         \vdots            & \vdots         & \ddots     & \vdots\\
 u_{n+1-i}^{n-2} &u_{n+1-i}^{n-3}u_{n-i}^{2}&\cdots &u_{n-i}^{n-2}\\
   \end{array} \right ] .
    \eqno(C4)
$$}

\bigskip

 \noindent
 Similar to $(A8)$, we consider the matrix $YQ_{n-1}^{-i}$. It is an easy matter to write
  down the matrix as multiplication by $Q_{n-1}$ shifts the indices (see $(A5)$ and $(A8)$).
  To calculate the determinant of $Y$, we
 consider the cofactor expansion of  the matrix $ A_1Q_n^{-i}$
 by the $i$th row (see $(A8)$),
  where the  $i$-th row of
   $A_1 Q_n^{-i}$ takes the form $( 1^{n-1},0, 0,\cdots, 0)$. An easy observation of the
    actual forms of the matrices $A_1Q_n^{-i}$ and $YQ_{n-1}^{-i}$ shows that

   $$\mbox{det}\, A_1Q_n^{-i} =(-1)^{i+1}(u_{1-i}u_{2-i}\cdots u_{-1}) (u_1 u_2\cdots u_{n-i})
   \mbox{det}\, YQ_{n-1}^{-i}.\eqno(C5)$$

   \medskip
   \noindent In summary,
    $$ C_i = \frac{
    W_{r+i+1}^n  \sigma(n)\sigma(n-1)\mbox{det}\,A_1\mbox{det}\,Q_{n-1}^i}
  {u_{i+1}u_i (u_{1-i}u_{2-i}\cdots u_{-1}) (u_1 u_2\cdots u_{n-i})
  \mbox{det}\,Q_n^{i}}.\eqno(C6)$$

 \medskip
 \noindent Recall that $u_{-m} = -q^{-m}u_m$ (this simple fact can be proved by induction)
  and that det$\,Q_n = q^{n(n-1)/2}$.
  This gives the determinant of $C$ (see $(C2)$).

\section {Appendix $D$}   Suppose that $u_1 u_2 \cdots u_v \ne 0$.
This appendix is devoted to the study of the characteristic polynomial of the matrix $
Q_{v}$ (see $(A4)$). For simplicity, we shall denote $Q_{v}$ by $Q$.
Let $E$ be the rational canonical form of $
 $$f(x) = \sum_{i=0}^{v} a_ix^i =  \sum_{i=0}^{v}
  (-1)^{i} q^{i(i-1)/2} (v|i)_u x^i.$ Then

  $$E=\left [
\begin{array}{ccccccc}
 0  &1& 0 & 0&\cdots  &0&0\\
 0&0&1&0&\cdots &0&0\\
0&0&0&1&\cdots &0&0\\
\vdots& \vdots& \vdots& \vdots&\ddots&\vdots&\vdots\\
0&0&0&0&\cdots  &0&1\\
 a_0&a_1& a_2& a_3&\cdots &a_{v-2}&a_{v-1}\\
   \end{array} \right ] .
    \eqno(D1)
$$
\medskip

\noindent
Let $A_1$ (set $n=v$) be the matrix given as in $(A2)$.
 Applying Proposition 3.2, the function $X(m) = u_{m-1}^k u_{m}^{v-1-k}$
 satisfies the recurrence relation
  $$
  \sum_{i=0}^{v}
  (-1)^{i} q^{i(i-1)/2} (v|i)_u
 u_{i-1}^k u_{i}^{v-1-k}
  =0.\eqno(D2)$$
\noindent
 Applying $(D2)$, the multiplication of $A_1$ by $E$ (to the left of $A_1$) shifts the indices.
  Since the multiplication of $A_1$ by $Q$ (to the right of $A_1$) shifts the indices as well
   (see $(A5)$), one has

 $$EA_1 = A_1Q ,\eqno (D3)$$

 \medskip
 \noindent
 which can be verified by direction calculation.
 Since the determinant of $A_1$ is nonzero (see $(A7)$), $E$ and $Q$ are similar to each other.
  As a consequence, the characteristic polynomial of $Q$ is $f(x)$.




\bigskip

\noindent MSC2010 : 11B39, 11B83.



\noindent




 \bigskip






\begin{thebibliography}{99}

\bibitem[B]{sloane}
T. A. Brenann, \emph{Fibonacci powers and Pascal's triangle in a matrix},
The Fibonacci
Quarterly \textbf{2} (1964), 93--103,117--184.


\bibitem[CK]{sloane}
C. Cooper  and R. E. Kennedy, \emph{Proof of a Result of Jarden by Generalizing a Proof of Carlitz},
The Fibonacci
Quarterly \textbf{33.4} (1995), 304--311.














\bibitem[H]{sloane}
A. F. Horadam, \emph{Basic Properties of a Certain Generalized Sequence of Numbers}, The Fibonacci
Quarterly \textbf{3.3} (1965), 161--76.





\bibitem[Ho]{sloane}
F. T. Howard, \emph{The sum of the squares of two generalized Fibonacci numbers}, The Fibonacci
Quarterly \textbf{41.1} (2003), 80--84.



\bibitem[J]{sloane}
D. Jarden, \emph{Recurring sequences}, Riveon Lematematika, (1958), 42--44.









\bibitem[LL1]{sloane}
C.L. Lang and M.L. Lang, \emph{Fibonacci Numbers and Identities}, to appear in
 The  Fibonacci Quarterly,
 {\tt arXiv:math/1303.5162v2 [math.NT]} (2013)


\bibitem[LL2]{sloane}
C.L. Lang and M.L. Lang, \emph{Fibonacci Numbers and Identities $II$}, preprint,
 {\tt arXiv:math/1304.3388v2 [math.NT]} (2013)











\bibitem[S]{sloane}
A. M. Stinchcombe, \emph{Recurrence relations for powers of recursion sequences},
  The Fibonacci
Quarterly \textbf{36.5} (1998), 443-447.






\bibitem[W]{sloane} E. W.  Weisstein,
{\em Fibinacci numbers,}
{\tt http:// www.mathworld.wolfram.com/FibonacciNumbers.html}.



\end{thebibliography}
\end{document}